\documentclass[letterpaper, 9pt, twocolumn]{extarticle}

\usepackage[utf8]{inputenc}

\usepackage[left=16mm, right=16mm, top=19mm, bottom=25mm]{geometry}
\usepackage{cite}
\usepackage{amsmath,amssymb,amsfonts}
\usepackage{algorithmic}
\usepackage{graphicx}
\usepackage{textcomp}
\usepackage{url}
\usepackage[normalem]{ulem}

\usepackage{color, xcolor} 
\usepackage{epsfig} 
\usepackage{amsmath, bbm, bm} 
\usepackage{amsfonts, amssymb}
\usepackage{mathrsfs}
\usepackage{anyfontsize}  

\usepackage{tikz}
\usetikzlibrary{intersections}
\usepackage{pgfplots}
\usepackage{setspace}

\makeatletter
\renewcommand\section{\@startsection {section}{1}{\z@}%
                                   {-1.5ex \@plus -1ex \@minus -.2ex}%
                                   {1.5ex \@plus .2ex}%
                                   {\normalfont\small\bfseries}}

\renewcommand\subsection{\@startsection{subsection}{2}{\z@}%
                                    {-3.25ex\@plus -1ex \@minus -.2ex}%
                                    {1.5ex \@plus .2ex}%
                                    {\normalfont\normalsize\itshape}}

\def\footnotemark{}

\renewcommand{\abstract}{
    \normalfont\small%
      \bfseries\textit{Abstract}---%
      \bfseries\selectfont
}

\makeatother

\newcommand {\st}{\text{\footnotesize $(t)$}}

\newcommand {\sta}[1]{\text{\footnotesize $(#1)$}}
\newcommand {\tc}{\tilde{c}}
\newcommand {\tI}{\tilde{I}}
\newcommand {\tR}{\tilde{R}}
\newcommand {\cB}{{\mathcal{B}}}
\newcommand {\tB}{\tilde{\cB}}
\newcommand {\vbeta}{\vec{\beta}}

\newcommand {\sS}{{\mathscr{S}}}
\newcommand {\sM}{{\mathscr{M}}}

\newcommand {\cL}{{\mathcal L}}

\newcommand {\cP}{{\mathcal{P}}}

\newcommand {\cS}{{\mathcal{S}}}
\newcommand {\cT}{{\mathcal{T}}}
\newcommand {\cV}{{\mathcal{V}}}

\newcommand {\cY}{{\mathscr{Y}}}

\newcommand {\R} {{\rm I\kern-2.5pt R}}
\newcommand {\C} {{\rm I\kern-5pt C}}

\newtheorem{assumption}{Assumption}

\newtheorem{theorem}{Theorem}
\newtheorem{remark}{Remark}
\newtheorem{example}{Example}
\newtheorem{definition}{Definition}

\newcommand{\beqa}{\begin{eqnarray}}
\newcommand{\eeqa}{\end{eqnarray}}
\newcommand{\beqan}{\begin{eqnarray*}}
\newcommand{\eeqan}{\end{eqnarray*}}
\newcommand{\beq}{\begin{equation}}
\newcommand{\eeq}{\end{equation}}

\newcommand{\bfl}{\begin{flushleft}}
\newcommand{\efl}{\end{flushleft}}
\newcommand{\bgnv}{\begin{verbatim}}
\newcommand{\ednv}{\end{verbatim}}

\newcommand{\myarr}{\begin{array}{lll}}

\newcommand{\bitem}{\begin{itemize}}
\newcommand{\eitem}{\end{itemize}}
\newcommand{\benum}{\begin{enumerate}}
\newcommand{\eenum}{\end{enumerate}}

\newcommand{\orcid}[1]{\href{https://orcid.org/#1}{\textcolor[HTML]{A6CE39}{\aiOrcid}}}

\newcommand{\la}[1]{{\bf {\color{red}[\marginpar[\hbox{{\large%
$\circ$}\raisebox{0ex}{\large $\longrightarrow$}}]%
{\hbox{\raisebox{0ex}{\large $\longleftarrow$}{\large
$\circ$}}}RL: #1]} } }

\newcommand{\Nuno}[1]{{\bf {\color{teal}[\marginpar[\hbox{{\large%
$\circ$}\raisebox{0ex}{\large $\longrightarrow$}}]%
{\hbox{\raisebox{0ex}{\large $\longleftarrow$}{\large
$\circ$}}}#1]} } }

\newcommand{\Jair}[1]{{\bf {\color{purple}[\marginpar[\hbox{{\large%
$\circ$}\raisebox{0ex}{\large $\longrightarrow$}}]%
{\hbox{\raisebox{0ex}{\large $\longleftarrow$}{\large
$\circ$}}}J: #1]} } }

\title{\Large \bf Epidemic Population Games With Nonnegligible Disease Death Rate} 

\author{Jair Cert\'{o}rio$^1$ and Nuno C. Martins$^1$ and Richard J. La$^1$
\thanks{Work supported by AFOSR (FA9550-19-1-0315) and NSF (2135561).}
\thanks{$^{1}$Electrical and Computer Engineering Department and Institute for Systems Research, University of Maryland at College Park.
        Email: {\tt\small \{certorio,nmartins,hyongla\}@umd.edu}}%
}

\pgfplotsset{compat=1.7}
\begin{document}

\maketitle
\thispagestyle{empty}
\pagestyle{empty}

\setlength{\voffset}{0.1 in}

\begin{spacing}{0.8}
\begin{abstract} 
    A recent article that combines normalized epidemic compartmental models and population games put forth a system theoretic approach to capture the coupling between a population's strategic behavior and the course of an epidemic. It introduced a payoff mechanism that governs the population's strategic choices via incentives, leading to the lowest endemic proportion of infectious individuals subject to cost constraints. Under the assumption that the disease death rate is approximately zero, it uses a Lyapunov function to prove convergence and formulate a quasi-convex program to compute an upper bound for the peak size of the population's infectious fraction. In this article, we generalize these results to the case in which the disease death rate is nonnegligible. This generalization brings on additional coupling terms in the normalized compartmental model, leading to a more intricate Lyapunov function and payoff mechanism. Moreover, the associated upper bound can no longer be determined exactly, but it can be computed with arbitrary accuracy by solving a set of convex programs.
\end{abstract}
\end{spacing}

\section{INTRODUCTION}  
    \label{sec:Introduction}

By fusing ideas from population games and epidemiological compartmental models, a recent article~\cite{Martins:2022um} proposed a system theoretic framework to capture the coupling between a population's strategic choices and the course of an epidemic.  There, as here, the population's agents can choose from $n$ strategies $\{1,\ldots,n\}$, each affecting differently the transmission rate of a susceptible-infectious-recovered-susceptible
(SIRS) compartmental epidemic model \cite{Hethcote00}. Each agent follows one strategy at a time, which it can repeatedly revise. A payoff vector in $\mathbb{R}^n$ structured as follows: 
\begin{equation} 
p\st:=r\st-c
\end{equation} quantifies the incentive (reward) $r\st$ provided by a planner for each strategy at time $t$, after the strategies' intrinsic costs $c$ are deducted. Namely, $c$ is the vector whose $\ell$-th entry $c_\ell$ is the inherent cost of the $\ell$-th strategy, and $r\st$ is a reward vector meant to incentivize the adoption of safer (costlier) strategies, where $r_\ell\st$ is $\ell$-th strategy's reward. 

\subsection{Brief Overview of~\cite{Martins:2022um}}

The main tenets in~\cite{Martins:2022um} are: {\bf(1)}~An evolutionary dynamics model (EDM) captures the agents' preferences in the revision process by specifying as a function of $p\st$ the rate at which strategies are selected or abandoned at time~$t$~(see~\S\ref{sec:EDM}). {\bf(2)}~A so-called epidemic population game (EPG)~(see~\S\ref{sec:EPG}) comprises {\bf(i)}~the SIRS model and {\bf(ii)}~a dynamic payoff mechanism that determines $r\st$. The ``inputs" of the dynamic payoff mechanism are the state of the SIRS model (epidemic variables) and the population state $x\st$ whose entries are the strategies' prevalence in the population at time $t$.

Using $\delta$-passivity concepts~\cite{Fox2013Population-Game}, in~\cite{Martins:2022um} the authors designed a payoff mechanism for which the population's infection prevalence converges globally asymptotically to its minimum, subject to cost constraints. The convergence is established via a Lyapunov function that is also used to determine upper-bounds for the peak size of the population's infectious fraction. As is explained in~\cite{Martins:2022um}, the exact knowledge of the (EDM) is not required for the payoff mechanism design nor is it needed to for the validity of the aforementioned convergence result and upper bound. Specifically, the design, convergence result and upper bound hold for any (EDM) which, via structural analysis, can be shown to be $\delta$-passive.

\subsection{Contributions and Comparison to~\cite{Martins:2022um}}
In this paper, we describe several modifications to the models and the payoff mechanism used in~\cite{Martins:2022um} that will allow us to consider a \underline{positive disease death rate} $\delta>0$, as opposed to assuming that it is negligible as was done in~\cite{Martins:2022um}.

Since the framework is rather similar, we will be able to reuse (via referencing) significant portions of~\cite{Martins:2022um}. Mostly, we will reproduce here the assumptions, definitions, and introductory portions strictly necessary to explain our work. 

The main technical differences between this article and~\cite{Martins:2022um} are: {\bf(1)}~Since, as in~\cite{Martins:2022um}, our SIRS model is normalized (by the population's size), the presence of nonnegligible disease death rate causes additional coupling terms and it perturbs symmetries. {\bf(2)}~These alterations make for a more intricate payoff mechanism and Lyapunov function. {\bf(3)}~While in~\cite{Martins:2022um} an upper bound for the peak size of the population's infectious portion is computed using a single-shot quasi-convex program, here an analogous bound can be determined only approximately via the solution of a set of convex programs. The approximation's accuracy can be made arbitrarily fine by increasing the number of convex programs~solved.

Notwithstanding these differences, the proof of our main convergence result~(Theorem~1) would mirror that for~\cite[Theorem~1]{Martins:2022um}. Consequently, instead of providing a complete (and somewhat redundant) proof for Theorem~1, in~\S\ref{sec:mainresult} we present a proof outline explaining how the main steps of the proof in~\cite{Martins:2022um} could be used here after appropriate changes.

\subsection{Related Literature} 
\label{subsec:Related}

Many studies in epidemiology make use of compartmental models. For example, Kermack and McKendrick used a deterministic model for modeling the transmissions in a closed population and showed the existence of a critical threshold density of susceptible individuals for the occurrence of a major epidemic. Today, this model is known as the susceptible-infectious-recovered (SIR) model. Other similar compartmental models introduce additional states, e.g.,  deceased (D), exposed (E), maternally-derived (M), vaccinated (V), and include SEIR/S, SIRD, SIRS, SIRV, and MSIR. We refer a reader to \cite{Anderson1991} for a comprehensive survey.

In closely related studies, di Lauro et al.~\cite{diLauro2021} and Sontag~\cite{Sontag2021An-explicit-for} studied the problem of identifying when non-pharmaceutical interventions (NPIs), e.g., quarantine and lockdowns, should be put in place to minimize the peak infections; \cite{diLauro2021} studied the optimal timing for one-shot intervention, whereas \cite{Sontag2021An-explicit-for} considered a fixed number of complete lockdowns. Al-Radhawi et al. \cite{Al-Radhawi2021} modeled 
NPIs during a prolonged epidemic as feedback effects and examined the problem of tuning the NPIs to regulate infection rates as an adaptive control problem. Using a singular-perturbation approach, the authors investigated the stability of disease-free and endemic steady states. Godara et al. \cite{Godara2021A-control-theor} considered the problem of controlling the infection rate to minimize the total cost till herd immunity is achieved in an SIR model. They formulated it as an optimal control problem subject to a constraint on the fraction of infectious population.

Even though we do not consider epidemic processes on general networks, their dynamics on networks are studied extensively (see \cite{Pastor-Satorras2015} for a review of the literature), including time-varying networks~\cite{Pare-TNSE18a}. Recently, the topic of mitigating disease or infection spread in a network enjoyed much attention, and researchers investigated optimal strategies using vaccines/immunization (prevention) \cite{Cohen-PRL03, Preciado-CDC13}, antidotes or curing rates (recovery) \cite{Borgs-RSA10, mai2018distributed, Ottaviano-JCN18} or a combination of both preventive and recovery measures \cite{Nowzari-TCNS17, Preciado-TCNS14}. 

A crucial aspect of controlling epidemic processes is human behavior and the strategic interactions among individuals, which determine their decisions over time in response to their (perceived) payoffs.
Game theory provides a natural framework and tools for studying such strategic interactions, and several recent studies adopted an {\em evolutionary} or {\em population game} framework \cite{Amaral2021, Arefin2020, Bauch2003, Bauch2004, Bauch2005, Chang2019, Cho2020, d'Onofrio2011, Kabir2019b, Kabir2020, Kuga2018, Wang2020}. We refer an interested reader to \cite{Chang2020} and references therein for a comprehensive survey of earlier studies.

Amaral et al.~\cite{Amaral2021} studied the effects
of perceived risks
when individuals can choose to 
voluntarily quarantine or continue their normal life. 
They showed that increased perceived risks result in
multiple infection peaks due to strategic interactions.
Kabir and Tanimoto~\cite{Kabir2020} considered
a similar setting and showed that naturally
acquired shield immunity is unlikely to be 
effective in suppressing an epidemic without 
additional social measures with low costs for
individuals.


\section{EVOLUTIONARY DYNAMICS MODEL  (EDM)}
\label{sec:EDM}

 Here, $x\st$ is the so-called {\it population state} taking values in $\mathbb{X}$ defined below and whose $\ell$-th entry $x_\ell\st$ is the proportion of the population adopting the $\ell$-th strategy at time $t$. 
\begin{equation*}
\mathbb{X}:= \Bigg \{x \in [0,1]^n \ \Big | \ \sum_{i=1}^n x_i =1 \Bigg \}
\end{equation*} 

Following the standard approach in~\cite[Section~4.1.2]{Sandholm2010Population-Game}, the following {\it evolutionary dynamics model} {\bf (EDM)} governs the dynamics of~$x$ in the large-population limit:
\begin{subequations}
\begin{equation}\tag{EDMa} \label{eq:EDM-DEF} \dot {x}\st = \mathcal{V} ( x\st,p\st ), \quad t\geq 0, 
\end{equation}
where the $i$-th component of $\mathcal{V}$ is specified as: 
\begin{multline} \tag{EDMb} \label{eq:EDMfromProtocol} \mathcal{V}_i ( x\st,p\st ) := \underbrace{\sum_{j=1, j \neq i}^{n}  x_j\st \mathcal{T}_{ji} ( x\st,p\st )}_{\text{\footnotesize inflow switching to strategy $i$}}  \\ - \underbrace{\sum_{j=1,j \neq i}^{n} x_i\st \mathcal{T}_{i j} ( x\st,p\st ) }_{\text{\footnotesize outflow switching away from strategy $i$}}
\end{multline}
\end{subequations}

A Lipschitz continuous map $\mathcal{T}: \mathbb{X} \times \mathbb{R}^{n} \rightarrow [ 0,\bar{\mathcal{T}}]^{n \times n}$, with upper bound $\bar{\mathcal{T}} > 0$, is referred to as {\it the revision protocol} and models the agents' strategy revision preferences. In \cite[Part~II]{Sandholm2010Population-Game} and \cite[\S 13.3-13.5]{Sandholm2015Handbook-of-gam} there is a comprehensive discussion on protocols types and the classes of bounded rationality rules they model. The analysis in \cite[\S{IV}]{Park2019From-Population} substantiates using (EDM) as a deterministic approximation for the case when a dynamical payoff mechanism generates $p$ from $x$, as will be the case~here. 

Below, we define a widely-used class of protocols, which we will repeatedly invoke to illustrate key concepts and employ in examples throughout the article.

\begin{definition} \label{def:IPC} Any protocol is said to be of the {\it impartial pairwise comparison} (IPC) type~\cite{Sandholm2010Pairwise-compar} if there is a map ${\phi:\mathbb{R}_{\geq 0} \rightarrow [0,\bar{\cT}]^n}$, whose components satisfy $\phi_j(0)=0$ and $\phi_j(\nu)>0$ for $\nu>0$, such that $\cT$ can be recast as:
\begin{equation}
\label{eq:defIPC}
\mathcal{T}_{ij}(x,p) \underset{\text{\tiny IPC}}{=} \phi_j([\tilde{p}_{ij}]_+),
\end{equation} where $\tilde{p}_{ij}:=p_j-p_i$. The so-called Smith's protocol~\cite{Smith1984The-stability-o} is ${\phi_j^{\text{\tiny Smith}}([\tilde{p}_{ij}]_+) := \min \{ \lambda [\tilde{p}_{ij}]_+, \bar{\mathcal{T}} \} }$, with $\lambda>0$.
\end{definition}

\subsection{Nash Stationarity and $\delta$-Passivity Assumption}
\label{subsec:DeltaPass}

We proceed to describe two assumptions on the (EDM) that will play a key role throughout the paper.

\begin{assumption} \label{ass:NS} {\bf (Nash Stationarity)}
We assume that $\cT$ is ``Nash~stationary", in the sense that the following holds:
\begin{equation}\tag{NS}
\cV(x,p)=0 \quad \Leftrightarrow \quad x \in \sM(p), \quad p \in \mathbb{R}^n
\end{equation} where $\sM : \mathbb{R}^n \rightarrow 2^\mathbb{X}$ is the following best response map:
\begin{equation*}
\sM(p) : = \arg \max_{x \in \mathbb{X} } \ p'x, \quad p \in \mathbb{R}^n. \end{equation*}
\end{assumption}

Therefore, for a revision protocol satisfying (NS), $x$ is an equilibrium of~(EDM) if and only if $x$ is a best response to~$p$. Any IPC protocol (see~Definition~\ref{def:IPC}) satisfies~(NS), as do other large classes of protocols~(see~\cite[\S13.5.3]{Sandholm2015Handbook-of-gam}).


Our analysis of the long-term evolution of the epidemic state and $(x,p)\st$ will leverage the following assumption stemming from the $\delta$-passivity concept originally proposed in~\cite{Fox2013Population-Game} and later generalized in~\cite{Park2019From-Population,Arcak2020Dissipativity-T,Kara2021Pairwise-Compar-a}.

\begin{assumption}\label{assm:deltaD} 
There exist a differentiable function ${\cS: \mathbb{X}\times \mathbb{R}^n \rightarrow \mathbb{R}_{\ge 0}}$ and a Lipschitz continuous function ${\cP: \mathbb{X}\times \mathbb{R}^n \rightarrow \mathbb{R}_{\ge 0}}$ that satisfy the following inequality for all $x$, $p$ and $u$ in $\mathbb{X}$, $\mathbb{R}^n$ and $\mathbb{R}^n$, respectively:
\begin{subequations}
\label{eq:delta-passivityWconditions}
\begin{equation}\label{delta-passivity}
\frac{\partial \cS(x,p)}{\partial x}\cV(x,p)+\frac{\partial \cS(x,p)}{\partial p}u \le -\cP(x,p)+u' \cV(x,p) 
\end{equation}
where $\cS$ and $\cP$ must also satisfy the equivalences below:
\begin{align}
\label{informative}
\cS(x,p)=0 \quad &\Leftrightarrow \quad \cV(x,p)=0. \\
\label{negdef}
\cP(x,p)=0 \quad &\Leftrightarrow \quad \cV(x,p)=0 
\end{align}

In addition, the following inequality (not required in standard $\delta$-passivity) must hold:
\begin{equation}
\label{homogeneity}
\cP(x,\alpha p) \geq \cP(x,p), \quad \alpha \geq 1, x \in \mathbb{X}, p \in \mathbb{R}^n
\end{equation}
\end{subequations}
\end{assumption}

Based on the Lyapunov functions in~\cite{Hofbauer2009Stable-games-an}, the authors of~\cite{Fox2013Population-Game,Kara2021Pairwise-Compar-a} determined, for main classes of protocols, explicit expressions for $\cS$ and $\cP$. The following is the main conclusion in \cite[Remark~4]{Martins:2022um}.

\begin{subequations}
\begin{remark} \label{rem:StructureForIPC}  Any IPC protocol~(\ref{eq:defIPC}) with non-decreasing $\{\phi_1, \ldots,\phi_n\}$ satisfies Assumption~\ref{assm:deltaD}.
\end{remark} 
\end{subequations}

\section{EPG WITH POSITIVE DISEASE DEATH RATE}
\label{sec:EPG}

Rather than focusing on what each strategy may represent, in our analysis we assume that a vector  $\vbeta$ in $\mathbb{R}_{>0}^n$ is given whose $\ell$-th entry $\vbeta_\ell$ quantifies the effect of strategy $\ell$ towards the transmission rate $\cB\st$ defined below:
\begin{equation}
\label{eq:betaaverage}
\cB\st := \vbeta' x\st, \quad t \geq 0
\end{equation}

The strategies' inherent costs decrease for higher transmission rates, and we order the entries of $\vbeta$ and $c$~as:
\begin{equation*}
 \vbeta_i < \vbeta_{i+1} \text{ and } c_i > c_{i+1}, \quad 1 \leq i \leq n-1
\end{equation*}
We will also use $\tilde{c}$ defined below:
\begin{equation*}
    \tilde{c}_i:=c_i-c_n, \quad 1 \leq i \leq n
\end{equation*}

\subsection{A Modified SIRS Model with Positive Disease Death Rate}
We proceed to describe a modification of the normalized SIRS model in~\cite{Martins:2022um} so as to allow for positive disease death rate $\delta>0$: let $\theta$ and $\zeta$ be the birth rate and the natural death rate (from all epidemic-unrelated causes), respectively. Similarly, $\gamma$ is the disease recovery rate, and $\psi$ is the rate at which a recovered individual becomes susceptible again due to waning immunity.

Suppose $N\st$ is the population size at time $t$, and let $I\st$, $R\st$ and $S\st$ be the numbers of infectious, recovered and susceptible individuals, respectively, at time $t$ divided by $N\st$. The population size $N\st$ is obtained as the solution of  $\dot{N}\st=(g-\delta I\st)N\st$, where $g := \theta - \zeta$. The epidemic model is then given by $I(t)$ and $R(t)$ (which are normalized by the population size $N(t)$)~\cite{Hethcote00}, where 
\begin{subequations} \label{eq:SIRS}
\begin{align}
    \dot{I}\st &= (\cB\st S\st  +\delta I\st - \sigma) I\st, \label{eq:SIRSa} \\ \label{eq:SIRSb}
    \dot{R}\st &= \gamma I\st - \omega R\st +\delta R\st I\st,
\end{align} 
\end{subequations} 
where $\bar{\sigma} := \gamma + \zeta + \delta$, $\sigma := g + \bar{\sigma}$, $\bar{\omega} = \psi + \zeta$ and $\omega:=g + \bar{\omega}$.
Our time unit is \underline{one day}, and $\bar{\sigma}^{-1}$ is the mean infectious period (in days) for an affected individual (till recovery or death). We also assume that newborns are susceptible.

\begin{assumption} \label{assmp:inequalities}
We consider that $\delta$ is positive, but moderate enough so that $\delta<\omega$ 
and $\delta < \gamma$. 
We also consider that $\sigma < \vbeta_1$, i.e., any strategy with transmission rate less than or equal to $\sigma$ would be unfeasible or too onerous. 
\end{assumption}

Notice that, in comparison to the SIRS model in~\cite{Martins:2022um}, there is an additional term in~(\ref{eq:SIRSa}) and a cross term in~(\ref{eq:SIRSb}), which will complicate our analysis. The models would be identical if one were to set $\delta=0$.

\subsection{Assumptions and Defining $\beta^*$ and $x^*$}

\begin{assumption} \label{assm:AboutBeta}
The following must hold when~$n \geq 3$:
\begin{equation}
\label{eq:BetaAndCIneq}
\frac{c_i-c_{i+1}}{\vbeta_{i+1}-\vbeta_i} > \frac{c_{i+1}-c_{i+2}}{\vbeta_{i+2}-\vbeta_{i+1}}, \quad 1 \leq i \leq n-2
\end{equation}
\end{assumption}

Under Assumption~\ref{assm:AboutBeta}, as the transmission rate decreases it becomes costlier to reduce it further. 
\begin{remark} \label{rem:KKTBetaStar} Using this assumption and Karush-Kuhn-Tucker conditions, we infer that for any given budget $c^*$ in $[0,\tc_1]\diagdown \{\tc_{n},\ldots,\tc_1\}$, the following optimization has a unique solution $x^*$ from which we also define $\beta^*$:
\begin{equation}
    x^*: = \arg \min \big \{ \vbeta' x \ | \ \tc'x \leq c^*, \ x \in \mathbb{X} \big \}, \quad \beta^*:=\vbeta'x^*
\end{equation} For $\tc_{i^*+1} < c^* < \tc_{i^*}$, the unique optimal solution is ${x^*_{i^*} = (c^*-\tc_{i^*+1})/(\tc_{i^*}-\tc_{i^*+1})}$, $x^*_{i^*+1}=1-x^*_{i^*}$ and the other entries of $x^*$ are zero, leading to $\vbeta_{i^*} < \beta^* < \vbeta_{i^*+1}$.
\end{remark}

As defined above, $\beta^*$ is the minimal transmission rate achievable for the given cost budget $c^*$. Notice that we exclude $\{\tc_{n},\ldots,\tc_1\}$ from the possible $c^*$ choices. This is a simplification intended to avoid the complications present in the analysis of (Case~II) in~\cite{Martins:2022um}.

\subsection{Lyapunov Stability of~(\ref{eq:SIRS}) When $\cB\st=\beta^*$}
We proceed to analyze the case when $\cB\st = \beta^* > \delta$ and
$(I^*,R^*) \neq (0,0)$ is the endemic equilibrium satisfying:
\begin{subequations}
\label{eq:EndemicEquil}
\begin{align}
    \beta^*-\sigma &= \beta^*(I^*+R^*)-\delta I^* \\
    0 &= \gamma I^* - \omega R^* +\delta R^* I^*
\end{align}
\end{subequations}

The upcoming Remark~\ref{rem:EndemicEq} can be used to ascertain the existence and uniqueness of $(I^*,R^*)$. To establish asymptotic stability of $(I^*,R^*)$ for~(\ref{eq:SIRS}), we use~(\ref{eq:EndemicEquil}) to rewrite~(\ref{eq:SIRS}) as:
\begin{subequations}
\label{eq:Param1SIRS}
\begin{align}
    \dot{I}\st &= (\beta^*(R^*-R\st) + (\beta^*-\delta)(I^*-I\st)) I\st, \\
    \dot{R}\st &= (\gamma +\delta R^*) (I\st-I^*) - (\omega-\delta I\st) (R\st-R^*),
\end{align} 
\end{subequations}

Now, consider the following candidate Lyapunov function based on a modification of the elegant one in~\cite{ORegan2010Lyapunov-functi}:
\begin{equation}
\label{eq:VDEF}
    \cV(I,R) : = (I-I^*)+I^* \ln \frac{I^*}{I} + \frac{a}{2}(R-R^*)^2
\end{equation} 
where $a:=\frac{\beta^*}{\gamma+\delta R^*}$.

The derivative of $\cV(I\st,R\st)$ along trajectories is:
\begin{multline}
\label{eq:cV-derivative}
    \tfrac{d}{dt} \cV(I\st,R\st) = -(\beta^*-\delta)(I\st-I^*)^2  \\ - a (\omega-\delta I\st)(R-R^*)^2
\end{multline}

Hence, if $\delta < \omega$ then (\ref{eq:cV-derivative}) is negative definite, which implies that $\cV$ is indeed a Lyapunov function in the domain of interest $(0,1]\times [0,1]$.

\subsection{A Useful Parameterization and EPG}
Start by defining $(\hat{I}_\cB,\hat{R}_\cB)$, for $\cB$ in $[\vbeta_1,\vbeta_n]$, as the unique solution in $(0,1] \times [0,1]$ of:
\begin{subequations}
\label{eq:IandRRoots2}
\begin{align}
    \cB-\sigma &= \cB (\hat{I}_\cB+\hat{R}_\cB)-\delta \hat{I}_\cB \\
    0 &= \gamma \hat{I}_\cB - \omega \hat{R}_\cB +\delta \hat{R}_\cB \hat{I}_\cB
\end{align}
\end{subequations} Namely, by solving a simple quadratic equation we obtain:
\begin{equation}
\label{eq:IandRRoots3}
    \hat{I}_\cB := \frac{b_\cB -\sqrt{\Delta}}{2\delta (\cB-\delta)}, \quad \hat{R}_\cB := (1-\frac{\sigma}{\cB}) - (1-\frac{\delta}{\cB}) \hat{I}_\cB
\end{equation} where $b_\cB:=\gamma \cB+\omega(\cB-\delta)+\delta(\cB-\sigma)$ and the discriminant is ${\Delta:=b_\cB^2 -4\delta\omega(\cB-\delta)(\cB-\sigma)}$. Notice that the endemic equilibrium solving~(\ref{eq:EndemicEquil}) is $(I^*,R^*)=(\hat{I}_{\beta^*},\hat{R}_{\beta^*})$. 
\begin{remark}\label{rem:EndemicEq} In the Appendix, we show that, subject to Assumption~\ref{assmp:inequalities}, $(\hat{I}_{\cB},\hat{R}_{\cB})$ is indeed the unique solution of~(\ref{eq:IandRRoots2}) in $(0,1]\times[0,1]$ and $\hat{I}_{\cB}+\hat{R}_{\cB} \leq 1$. 
\end{remark}
In the same way we derived~(\ref{eq:Param1SIRS}), we can rewrite~(\ref{eq:SIRS}) as:
\begin{subequations}
\label{eq:EPG}
\begin{align}
    \dot{I}\st &= (\cB\st\tR\st + (\cB\st-\delta)\tI\st) I\st, \\
    \dot{R}\st &= (\omega-\delta I\st) \tR\st-(\gamma +\delta \hat{R}_{\cB(t)}) \tI\st,
\end{align} where $\tR\st:=\hat{R}_{\cB(t)}-R\st$ and $\tI\st:=\hat{I}_{\cB(t)}-I\st$.
 The following dynamic payoff mechanism is adapted from~\cite{Martins:2022um}:
\begin{align}
\dot{q}\st &= G(I\st,R\st,x\st,q\st) \\ r\st &= q\st \vbeta + r^*, \text{or equivalently,} \ p\st=q\st \vbeta+r^o,
\end{align} where $q\st \in \mathbb{R}$, $G$ is a Lipschitz continuous map we seek to design, $r^o:=r^*-c$, and $r^*$ is chosen to satisfy:
\begin{equation}    \label{eq:r*}
    \begin{cases} r^*_i < \tilde{c}_i & \text{if $n>2$ and $x_i^*=0$} \\
    r^*_i = \tilde{c}_i & \text{otherwise} \end{cases}, \quad 1 \leq i \leq n
\end{equation}
\end{subequations}

\begin{definition} {\bf Epidemic population game (EPG)} We refer to the dynamical system specified in~(\ref{eq:EPG}) as an \emph{ epidemic population game (EPG)}. The state of the EPG is ${Y:=(I,R,x,q)}$ and takes values in $\mathbb{Y}$ specified below:
\begin{equation}
    \mathbb{Y}:=\{(I,R) \ | \ 0<I\leq 1, 0\leq R\leq 1-I\} \times \mathbb{X} \times \mathbb{R}
\end{equation}
\end{definition}

\section{LYAPUNOV FUNCTION AND SPECIFYING $G$}
We start by defining a candidate Lyapunov function:
\begin{equation}    \label{eq:cL}
    \cL(Y) : = \sS(I,R,\cB)+\cS(x,p), \quad Y \in \mathbb{Y}
\end{equation} 
where $\cS$ satisfies Assumption~\ref{assm:deltaD}, $\cB:=\vbeta'x$, and $\sS$, which is a modification of $\cV$ in~(\ref{eq:VDEF}), is defined as follows:
\begin{equation}    \label{eq:sS}
    \sS(I,R,\cB):=  \hat{I}_\cB \ln \frac{\hat{I}_\cB}{I} -\tI+ \frac{a_\cB}{2}\tR^2+\frac{\upsilon^2}{2}\tB^2
\end{equation} with $a_\cB : = \frac{\cB}{\gamma+\delta \hat{R}_\cB}$ and $(\tI,\tR,\tB):=(\hat{I}_{\cB}-I,\hat{R}_{\cB}-R,\cB-\beta^*)$.

\begin{remark} Notice that $\sS(I,R,\cB) \geq 0$ for $Y$ in $\mathbb{Y}$ and $\sS(I,R,\cB)=0$ if and only if $(I,R,\cB)=(I^*,R^*,\beta^*)$. In addition, for any given $\cB$ in the interval $(\vbeta_1,\vbeta_n)$, the function $\sS(I,R,\cB)$ is convex with respect to $(I,R)$ in $(0,1]\times[0,1]$.
\end{remark}

\subsection{Using $\cL$ to Obtain a Stabilizing Policy $G$}
We now proceed to choose $G$ in such a way that the derivative of $\cL(Y\st)$ is non-positive along trajectories. Namely, we select $G(I,R,x,q)=-\partial_\cB \sS(I,R,\cB)$, where $\partial_\cB$ indicates the partial derivative with respect to $\cB$. After computing the partial derivatives, we obtain the following expression for $G$:
\begin{align} \label{eq:Gformula}
    G(I,R,x,q) &= \left ( \ln \frac{I}{\hat{I}_\cB} \right ) \partial_\cB\hat{I}_\cB-\upsilon^2\tB-G_\cB(R,\cB)\tR \\ \nonumber
    G_\cB(R,\cB) &: = \frac{1}{2}(2a_\cB\partial_\cB\hat{R}_\cB+(\hat{R}_\cB-R)\partial_\cB a_\cB)
\end{align} Furthermore, $\partial_\cB a_\cB = (\gamma+\delta \hat{R}_\cB)^{-2}(\gamma+\delta (\hat{R}_\cB-\cB\partial_\cB \hat{R}_\cB))$, and $(\partial_\cB\hat{I}_\cB,\partial_\cB\hat{R}_\cB)$ can be determined from~(\ref{eq:IandRRoots2}) as:
\begin{equation} \label{eq:partials}
    \begin{bmatrix} \partial_\cB\hat{I}_\cB \\ \partial_\cB\hat{R}_\cB \end{bmatrix} = \begin{bmatrix} \cB-\delta & \cB \\ \gamma+\delta \hat{R}_\cB & -(\omega - \delta \hat{I}_\cB) \end{bmatrix}^{-1} \begin{bmatrix} 1-\hat{I}_\cB-\hat{R}_\cB \\ 0 \end{bmatrix}
\end{equation}

\begin{remark} \label{rem:boundednessofpartials} We infer from Remark~\ref{rem:EndemicEq} that, subject to Assumption~\ref{assmp:inequalities}, the inverse matrix in~(\ref{eq:partials}) exists and
$( \partial_\cB\hat{I}_\cB, \partial_\cB\hat{R}_\cB)$ is uniformly bounded for all $\cB$ in~$[\vbeta_1,\vbeta_n]$. Furthermore, $G(I^*,R^*,x^*,q)=0$ for any $q$ in $\mathbb{R}$.
\end{remark}  

\subsection{Stability Notion and Main Result}
\label{sec:mainresult}

We adopt the following global asymptotic stability notion. 
\begin{definition} \label{def:GAS}  We say that $e^*$ in $\mathbb{Y}$ is globally asymptotically stable {\bf (GAS)} if it satisfies the following two conditions with respect to $\cL$ and solutions of the system formed by (EDM) and~(EPG): 
{\bf (i)}~It holds that ${\cY = e^*\Leftrightarrow \cL(\cY) = 0}$; {\bf (ii)}~For any $\cY\sta{0}$ in $\mathbb{Y}$, $e^*$ is the one and only accumulation point of $\{ \cY\st \ | \ t\geq 0\}$.
\end{definition}



The following theorem is our main result.

\begin{theorem} \label{thm:MainTheorem}
Let the protocol defining~(EDM) and the design parameters $\upsilon>0$ and $c^*$ in $(0,\tc_1) \diagdown \{\tc_{n},\ldots,\tc_1\}$ be given. If (NS) and Assumptions~\ref{ass:NS}-\ref{assmp:inequalities} hold, then for $G$ given by~(\ref{eq:Gformula}) the equilibrium $e^*:=(I^*,R^*,x^*,0)$ is GAS.
\end{theorem}
Hence, if the Theorem's conditions hold then:
\begin{align}
    \lim_{t \rightarrow \infty} (I,R,x,q,\cB)\st &= (I^*,R^*,x^*,0,\beta^*) \\
    \lim_{t \rightarrow \infty} r\st'x\st & = c^*
\end{align} 
where $t^{-1}\int_0^t r\sta{\tau} 'x\sta{\tau} d\tau$ quantifies the control policy's average cost up until $t$ and whose limit is also $c^*$.
Notice that the discussion about the universality of~\cite[Theorem~1]{Martins:2022um} remains valid here, and that, as was the case in~\cite{Martins:2022um}, we will be able to use $\cL$ in \S\ref{sec:example} to construct anytime bounds.

{\bf Outline of a Proof for Theorem~\ref{thm:MainTheorem}:}
Under the conditions of the theorem, and after some algebra, we obtain:
\begin{multline} \label{eq:derivativeLyapunov}
    \tfrac{d}{dt} \cL(Y\st) \leq -\cP(x\st,p\st)  -(\cB\st-\delta)\tilde{I}\st^2 \\  - a_{\cB\st} (\omega-\delta I\st)\tilde{R}\st^2
\end{multline}

Using~(\ref{eq:derivativeLyapunov}), the proof of the theorem can be carried out by following the same LaSalle-type method described in steps~1-4 in~\cite[Appendix~A (for Case~I)]{Martins:2022um}, which need not be replicated here. In following these steps to prove the theorem one would need to observe the following facts: {\bf (1)}~Although $\cL$ defined here is distinct from that in~\cite{Martins:2022um}, under Assumption~\ref{assmp:inequalities}, it can be used in the same way once one realizes that the second and third terms of the right hand side of (\ref{eq:derivativeLyapunov}) are zero only when $\tilde{R}=\tilde{I}=0$, and $\cP$ satisfies the same properties as in~\cite{Martins:2022um} (see~Assumption~\ref{assm:deltaD}).  {\bf (2)}~$G$ used here is also different from that in~\cite{Martins:2022um}, but $G$ would be used directly only when replicating step~2. To do so, we would need to use Remark~\ref{rem:boundednessofpartials} and carefully analyze $\hat{I}_\cB$, $a_\cB$ and $\partial_\cB a_\cB$ to conclude from~(\ref{eq:Gformula}) that, given any $\xi > 0$, there is $\mu_\xi>0$ such that the following holds for any $Y$ in $\mathbb{Y}$:
\begin{equation} \label{eq:ineqforproof}
  \tilde{I}^2+\tilde{R}^2 < \mu_\xi \implies  | G(I,R,x,q) + \upsilon^2\tilde{\cB}| < \tfrac{1}{3} \xi
\end{equation}
Such an analysis would be straightforward once one ascertained that  $a_\cB$, $\partial_\cB a_\cB$ and $\ln \hat{I}_\cB$, are uniformly bounded for $\cB$ in $[\vbeta_1,\vbeta_n]$ (see~Appendix). Once~(\ref{eq:ineqforproof}) is established one could complete step~2 to conclude that $q\st$ is bounded. $\square$ 

\section{USING $\cL$ TO OBTAIN AN ANYTIME BOUND}
\label{sec:example}

By~(\ref{eq:derivativeLyapunov}) and~(\ref{eq:cL}), the following bounds hold for all $t \geq 0$:
\begin{equation}
\label{eq:LyapunovBound}
\cL(Y\sta{0}) \geq \cL(Y\st) \geq \sS(I\st,R\st,\cB\st)
\end{equation}  
Recall that, although we omitted the dependence, $\sS$ depends on the design parameter $\upsilon$.

As mentioned in \S\ref{sec:Introduction}, in \cite{Martins:2022um} an upper bound on $I(t)$ could be computed from a quasi-convex program. Here, because $\sS$ is convex only in $(I, R)$ for fixed $\cB$, we cannot obtain an upper bound on $I\st$ by solving a single quasi-convex program; instead, an analogous bound can be determined only approximately via the solution of a set of convex programs as follows: 
for $\alpha > 0$, define
\begin{equation} \label{eq:piStarDef}
\pi^*_\upsilon(\alpha)
:= \sup \big\{ \ I / I^* \ | \ \sS(I, R, \vbeta' x)  \leq \alpha, \ Y \in \mathbb{Y} \ \big\}. 
\end{equation}

Consider the case ${n \! = \! 2}$ with an initial endemic equilibrium point $Y\sta{0}$ : $(I\sta{0}, R\sta{0})$ ${= \! (\hat{I}_{\cB\sta{0}}, \hat{R}_{\cB\sta{0}})}$ and ${q\sta{0} \! = \! 0}$, where ${\cB\sta{0} \! = \! \vbeta' x\sta{0}}$. Because ${p\sta{0} \! = \! r^o}$, both strategies have a payoff of $- c_n$ at time ${t \! = \! 0}$ and ${\dot{x}\sta{0} \! = \! \cV(x\sta{0}, p\sta{0}) \! = \! {\bf 0}}$, which implies ${\cS(x\sta{0}, p\sta{0}) \! = \! 0}$ from \eqref{informative}. Thus, from \eqref{eq:sS}, we have ${\sS(I\sta{0}, R\sta{0}, \cB\sta{0}) \! = \! 0.5 \upsilon^2 (\cB\sta{0} - \beta^*)^2}$. It now follows from the inequalities in \eqref{eq:LyapunovBound}
and the definition of $\pi^*_{\upsilon}(\alpha)$ in \eqref{eq:piStarDef} that 
\begin{equation} \label{eq:piStarBound}
I\st \leq I^* \text{\tiny $\times$} \pi^*_{\upsilon} \left( \tfrac{1}{2} \upsilon^2 (\cB\sta{0} - \beta^*)^2 \right ).
\end{equation}
Following a similar argument, an analogous bound on $I\st$ can be obtained for the case with $n \geq 3$.

In practice, computing $\pi^*_{\upsilon}(\alpha)$ exactly is problematic. But, it can be approximated as we explain here: $\pi^*_{\upsilon}(\alpha) = \max\{ \  \pi_\upsilon^{\cB}(\alpha) \ | \ \cB \in [\vbeta_1, \vbeta_n] \ \}$, where
\begin{equation*}    \label{eq:piB}
\pi^\cB_\upsilon(\alpha)
:= \sup \big\{ \ I / I^* \ | \ \sS(I, R, \cB)  \leq \alpha, \ \vbeta' x = \cB, \ Y \in \mathbb{Y} \ \big\}. 
\end{equation*}
Note that this optimization problem is convex and can be solved efficiently. Now, $\pi^*_{\upsilon}(\alpha)$ can be approximated using $\tilde{\pi}_{\upsilon}(\alpha) = \max\{ \  \pi_\upsilon^{\cB^{(m)}}(\alpha) : m = 1, \ldots, M \ \}$, where $\{\cB^{(m)} : m = 1, \ldots, M\}$ is a suitably chosen finite subset of $[\vbeta_1, \vbeta_n]$.

We will use the following example to illustrate the validity of our bounds. \underline{Our time unit will be one day.}

\begin{example} \label{example1} Consider an example with following parameters: $g=0$, $\gamma=0.1$ (mean recovery period $\sim$ 10 days), ${\sigma=\bar{\sigma}=0.105}$, $\delta = 0.005$, and $\omega = \bar{\omega}=0.011$~(mean immunity period $\sim$ 91 days). 
In addition, $(\vbeta_1, \vbeta_2)=(0.15, 0.19)$ and $(c_1, c_2)=(0.2, 0)$. We select $c^*=0.1$, which yields $x^*_1=x^*_2=0.5$, $\beta^*=0.17$, and $(I^*,R^*) \approx (6.89 \%,37.52 \%) $. We assume $x_1\sta{0}=1$, $\cB\sta{0}=0.15$, and $(I\sta{0},R\sta{0}) = (\hat{I}\sta{0},\hat{R}\sta{0})=(2.94\%,27.15\%)$. Our goal is to design $G$ and $H$ so that ${I\st} \leq 1.3 \times {I^*}$ for all $t \geq 0$. 
\end{example}

\begin{figure}
    \centering
    \vspace{.2in}
    \input{fig1_v3.tikz}
    \caption{Plot of $\tilde{\pi}_{\upsilon} \big( \tfrac{1}{2}\upsilon^2(\cB\sta{0} - \beta^*)^2 \big)$ in Example~1 as a function of $\upsilon$ for varied $\beta^*$ and $\delta$   (other parameters of Example~1 are unchanged.) With $M=30$ and equidistant points in $[\vbeta_1, \vbeta_n]$ to approximate $\pi_\upsilon^*$.}
    \label{fig:pistar}
\end{figure}

\begin{figure}
    \centering
    \input{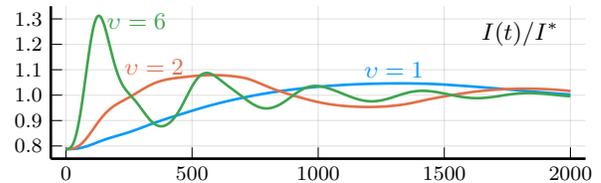}
    \caption{Simulation for Example~\ref{example1} using 
    $\upsilon$ as shown, and a Smith's protocol specified by $\lambda=0.1$ and $\bar{\cT}=0.1$. } \label{fig:simulationI}
\end{figure}

Analogously to the example in \cite{Martins:2022um}, this example shows the effect of reducing the budget from $r\sta{0}' x\sta{0} = 0.2$ to $c^* = 0.1$, given the initial conditions described above.
From Fig.~\ref{fig:pistar}, for $\beta^* = 0.17$ and $\delta = 0.005$, we have $\tilde{\pi}_2(\frac{1}{2} \cdot 2^2 \cdot 0.02^2) \approx 1.3$. Thus, by \eqref{eq:piStarBound}, 
when $\upsilon$ is chosen to be smaller than 2, we will be able to satisfy the requirement $I\st \leq 1.3 \times I^*$.
%
Although smaller values of $\upsilon$ reduce the overshoot, Fig.~\ref{fig:simulationI} shows that it also leads to slower convergence, which may keep the instantaneous cost, $r\st' x\st$,  larger for a longer period. Hence, in order to bring down the instantaneous cost quickly, $\upsilon$ should be chosen so that $\pi^*_{\upsilon}\big( \frac{1}{2} \upsilon^2(\cB\sta{0} - \beta^*)^2 \big)$ is close to the target overshoot. 

Fig.~\ref{fig:simulationI} also suggests that the bound in \eqref{eq:piStarBound} may be conservative for Smith's protocol. This is likely due to the fact that our bound is required to hold for {\em any} protocol satisfying the conditions in Theorem~\ref{thm:MainTheorem}. We note that choosing $\upsilon = 6$ no longer satisfies the target bound on $I\sta{t}$ (of $1.3 \times I^*$). 

\section*{APPENDIX}

Below, we explain why~(\ref{eq:IandRRoots2}) has a unique solution in $(0,1]\times[0,1]$ and why $\ln \hat{I}_\cB$ is uniformly bounded.

Using completion of squares and the facts that $\delta <\cB$ and $\sigma < \cB $, we infer that $\Delta$ in~(\ref{eq:IandRRoots3}) satisfies the following inequality, which also implies that $\hat{I}_\cB$ is real and positive: 
\begin{equation}
    b_\cB > \Delta \geq (\gamma \cB)^2 +2\gamma \cB(\omega(\cB-\delta)+\delta(\cB-\sigma))
\end{equation}


A simple analysis would lead to the following equivalence:
\begin{equation}
       \left ( \tfrac{b_\cB-\sqrt{\Delta}}{2\delta (\cB-\delta)} \geq 1 \text{ or } \tfrac{b_\cB+\sqrt{\Delta}}{2\delta (\cB-\delta)} \! \leq \! 1 \right ) \Leftrightarrow \gamma \cB \leq (\delta-\omega) (\sigma-\delta)
\end{equation} 
implying that there is always a unique $\hat{I}_\cB \in (0,1)$ such that (\ref{eq:IandRRoots2}) holds.
%
%
If $\hat{R}_\cB<0$, by \eqref{eq:IandRRoots2} we have that $0 = \gamma \hat{I}_\cB + \hat{R}_\cB (\delta \hat{I}_\cB - \omega )$, which is a contradiction since $\omega > \delta$ and $ 0<\hat{I}_\cB \leq 1$. If $\hat{R}_\cB>1-\hat{I}_\cB$, by \eqref{eq:IandRRoots3} we have
$\delta \hat{I}_\cB -\sigma >0$, also a contradiction since $\sigma = \bar{\sigma} + g > \delta$. Hence $0 \leq \hat{R}_\cB \leq 1-\hat{I}_\cB$ and there will always be a unique solution for (\ref{eq:IandRRoots2}) in $[0,1]^2$. Finally, we can show from Assumption~\ref{assmp:inequalities} that $\ln(\hat{I}_{\cB}(t))$ is uniformly bounded because $\hat{I}_{\cB}(t) \geq \Delta^* / \big(2 \delta (\beta_n - \delta) \big) > 0$, where $\Delta^* = b_{\cB}^* - \sqrt{(b_{\cB}^*)^2 - 4 \delta \omega(\beta_1 - \delta)(\beta_1 - \sigma)}$, and $b_{\cB}^* = \gamma \beta_n + \omega(\beta_n - \delta) + \delta(\beta_n - \sigma)$.

\bibliographystyle{ieeetr}
\bibliography{MartinsRefs,La}

\end{document}